\documentclass[final]{siamltex}

\usepackage{amsmath}
\usepackage{amsfonts}
\usepackage{graphicx}
\usepackage{booktabs}
\usepackage{tabularx}
\usepackage{url}

\def\clap#1{\hbox to 0pt{\hss#1\hss}}

\graphicspath{{.}{./figs/}}

\newcommand{\mA}{\mathbf{A}}

\newcommand{\mI}{\mathbf{I}}
\newcommand{\mJ}{\mathbf{J}}
\newcommand{\mP}{\mathbf{P}}
\newcommand{\mQ}{\mathbf{Q}}

\newcommand{\mU}{\mathbf{U}}
\newcommand{\mV}{\mathbf{V}}
\newcommand{\mW}{\mathbf{W}}
\newcommand{\mX}{\mathbf{X}}
\newcommand{\mY}{\mathbf{Y}}
\newcommand{\mZ}{\mathbf{Z}}

\newcommand{\ve}{\mathbf{e}}

\newcommand{\vq}{\mathbf{q}}
\newcommand{\vr}{\mathbf{r}}
\newcommand{\vu}{\mathbf{u}}
\newcommand{\vv}{\mathbf{v}}
\newcommand{\vw}{\mathbf{w}}
\newcommand{\vx}{\mathbf{x}}
\newcommand{\vy}{\mathbf{y}}

\newcommand{\sI}{\mathcal{I}}
\newcommand{\sJ}{\mathcal{J}}
\newcommand{\sS}{\mathcal{S}}

\newcommand{\ip}[1]{\left\langle #1 \right\rangle}
\newcommand{\tvec}{\mathrm{vec}}

\newcommand{\blambda}{\mbox{\boldmath$\lambda$}}
\newcommand{\bpi}{\mbox{\boldmath$\pi$}}

\newcommand{\mytitlenote}{This paper appeared in SIAM Journal of Scientific Computing in 2010 under the title \textit{A factorization of the spectral Galerkin system for parameterized matrix equations: derivation and applications}~\cite{sisc-galerkin2010}. We have changed the title of the paper here on \url{arXiv.org} to increase its prominence in relevant internet search results. Beyond the title and this footnote, this paper is identical to the published manuscript.}
\title{FACTORIZING THE STOCHASTIC GALERKIN SYSTEM
\footnote{ \protect\mytitlenote }}

\author{Paul G.~Constantine\thanks{Sandia National Labs,
Albuquerque, New Mexico. ({\tt pconsta@sandia.gov}).}
\and 
David F.~Gleich\thanks{Sandia National Labs, Livermore, California. ({\tt dfgleic@sandia.gov}).}
\and
Gianluca Iaccarino\thanks{Stanford University, Stanford, California ({\tt jops@stanford.edu}).}
}

\begin{document}

\maketitle

\begin{abstract}
Recent work has explored solver strategies for the linear system of equations arising from a spectral Galerkin
approximation of the solution of PDEs with parameterized (or stochastic) inputs. We consider the related problem of a
matrix equation whose matrix and right hand side depend on a set of parameters (e.g.\ a PDE with stochastic inputs
semidiscretized in space) and examine the linear system arising from a similar Galerkin approximation of the solution.
We derive a useful factorization of this system of equations, which yields bounds on the eigenvalues, clues to
preconditioning, and a flexible implementation method for a wide array of problems. We complement this analysis with
(i) a numerical study of preconditioners on a standard elliptic PDE test problem and (ii) a fluids application using
existing CFD codes; the \textsc{MATLAB} codes used in the numerical studies are available online.
\end{abstract}


\pagestyle{myheadings}
\thispagestyle{plain}
\markboth{P.~G. CONSTANTINE, D.~F. GLEICH, AND G. IACCARINO}{STOCHASTIC GALERKIN FACTORIZATION}

\section{Introduction}
Complex engineering models are often described by systems of equations, where the model outputs depend on a set
of input parameters. Given values for the input parameters, the model output may be computed by solving a
discretized differential equation, which typically involves the solution (or series of solutions) of a system of linear
equations for the unknowns. Each of these computations may be very expensive depending on factors such as grid
resolution or physics components in the model. As the role of simulation gains prominence in areas such as
decision support, design optimization, and predictive science, understanding the effects of variability in the input
parameters on variability in the model output becomes more important. Exhaustive parameter studies (e.g.\ uncertainty
quantification, sensitivity analysis, model calibration) can be prohibitively expensive -- particularly for a large
number of input parameters -- and therefore accurate interpolation and robust surrogate models are essential.

We consider the model problem of a matrix equation whose matrix and right hand side depend on a set of
parameters. Let $s\in\sS$ be a set of input parameters from a $d$-dimensional tensor product parameter space
$\sS=\sS_1\otimes\cdots\otimes\sS_d$; the range of each $\sS_i$ may be bounded or unbounded. We equip the parameter
space with a bounded, separable, positive weight function $\omega:\sS \mapsto \mathbb{R}_+$, where
$\omega(s)=\omega_1(s_1)\cdots\omega_d(s_d)$. (In a probabilistic context, this weight function represents a
probability measure on the input parameter space.) Let $A(s)$ be an $N\times N$ matrix-valued function where we assume
that $A(s)$ is invertible for all $s\in\sS$, and let $b(s)$ be the vector-valued function with $N$ components, where each
component is square integrable with respect to $\omega$. We seek the vector valued function $x(s)$ that
satisfies
\begin{equation}
\label{eq:parmmat}
A(s)x(s)=b(s),\qquad s\in\sS.
\end{equation}
Such parameterized matrix problems often arise as an intermediate step when computing an approximate
solution of a complex model with multiple input parameters. They appear in such diverse fields
as differential equations with random (or parameterized) inputs~\cite{Babuska01,Prudhomme02}, electronic circuit
design~\cite{Li2009}, image deblurring models~\cite{Chung08}, and ranking methods for nodes in a
graph~\cite{Page99,Constantine10b}. Given a parameterized matrix equation, one may wish to approximate the vector
valued function that satisfies the equation or estimate its statistics. Once an approximation is constructed that is
cheaper to evaluate than the true solution, statistics -- such as mean and variance -- of the approximation represent
estimates of statistics of the true solution.
 
The approximation model is a vector of multivariate polynomials represented as a series of orthonormal
polynomial basis functions where each basis function is a product of univariate orthonormal polynomials. We employ the
standard multi-index notation; let $\alpha=(\alpha_1,\dots,\alpha_d)\in\mathbb{N}^d$ be a multi-index, and define
the basis polynomial
\begin{equation}
\pi_\alpha(s) = \pi_{\alpha_1}(s_1)\cdots\pi_{\alpha_d}(s_d).
\end{equation} 
The polynomial $\pi_{\alpha_i}(s_i)$ is the orthonormal polynomial of degree $\alpha_i$, where the orthogonality is
defined with respect to the weight function $\omega_i(s_i)$. Then for $\alpha,\beta\in\mathbb{N}^d$,
\begin{equation}
\int_{\sS} \pi_\alpha(s)\pi_\beta(s)\omega(s)\,ds = 
\left\{
\begin{array}{cl}
1,&\alpha=\beta\\
0,&\mbox{otherwise}
\end{array}
\right.
\end{equation}
where equality between multi-indicies means component-wise equality. For a given index set $\sI\subset\mathbb{N}^d$ with
size $|\sI|<\infty$, the polynomial approximation can be written
\begin{equation}
\label{eq:polymodel}
x(s)\;\approx\;\sum_{\alpha\in\sI} \vx_\alpha \pi_\alpha(s) \;=\; \mX\bpi(s).
\end{equation}
The $N$-vector $\vx_\alpha$ is the coefficient of the series corresponding to $\pi_\alpha(s)$. The $N\times|\sI|$ matrix
$\mX$ has columns $\vx_\alpha$, and the parameterized vector $\bpi(s)$ contains the basis polynomials. The goal of the
approximation method is to compute the unknown coefficients $\mX$.  

Such polynomial models have become popular for approximating the solution of PDEs with random
inputs~\cite{LeMaitre10,Xiu2010}; they appear under names such as polynomial chaos methods~\cite{Xiu02}, stochastic
finite element methods~\cite{Ghanem91}, stochastic Galerkin methods~\cite{Babuska04}, and stochastic collocation
methods~\cite{Xiu05,Babuska07}. The Galerkin methods compute the series coefficients such that the equation residual is
orthogonal to the approximation space defined by the index set $\sI$; they typically employ a full polynomial basis of
order $n$ given by
\begin{equation}
\sI\;=\; \sI_n \;=\; \{\alpha\in\mathbb{N}^d\;:\;\alpha_1+\cdots+\alpha_d\leq n\},
\end{equation}
although such a basis set is not strictly necessary. The number of terms in this basis set is
$|\sI_n|=\binom{n+d}{n}$, which grows rapidly for $d>1$. The process of computing
the coefficients involves solving a linear system of size $N|\sI|\times N|\sI|$, which can be prohibitively large for
even a moderate number of input parameters (6 to 10) and low order polynomials (degree $< 5$). For this reason, there
has been a flurry of recent work on solver strategies for the matrix equations arising from the Galerkin
methods~\cite{Pellissetti00,Rosseel08,Keese05,Elman07,Elman05,Ernst09,Rosseel10}, including papers on
preconditioning~\cite{Powell09-2,Elman10,Powell08,Ullmann2010}. Such work has relied on knowing the matrix valued
coefficients $\mA_\alpha$ of a series expansion of the parameterized matrix $A(s)$,
\begin{equation}
\label{eq:Aexpansion}
A(s) = \sum_\alpha \mA_\alpha \pi_\alpha(s),
\end{equation}
which typically come from a specific form of the coefficients of the elliptic operator in the PDE models. Another
drawback of the Galerkin method is its limited ability to take advantage of existing solvers for the problem
$A(\lambda)x(\lambda)=b(\lambda)$ given a parameter point $\lambda\in\sS$. In contrast, the pseudospectral and
collocation methods can compute coefficients of the model \eqref{eq:polymodel} using only evaluations of the solution
vector $x(\lambda)$. The chosen
points typically correspond to a multivariate quadrature rule, and computing the coefficients of the polynomial model
\eqref{eq:polymodel} becomes equivalent to approximating its Fourier coefficients with the quadrature
rule~\cite{Babuska07,Constantine10c}. Thus, from the point of view of code reuse and rapid implementation, the
pseudospectral/collocation methods have a distinct advantage.

We propose a variant of the Galerkin method that alleviates the drawbacks of limited code reuse and memory
limitations. By formally replacing the integration in the Galerkin method by a multivariate quadrature rule -- a
step that is more often than not performed in practice -- we derive a decomposition of the linear system of equations
used to compute the Galerkin coefficients of \eqref{eq:polymodel}; such a method has been called Galerkin with
numerical integration in the context of numerical methods for PDEs~\cite{Canuto06}. This decomposition allows us to
compute the Galerkin coefficients using only evaluations of the parameterized matrix $A(\lambda)$ and parameterized
vector $b(\lambda)$ for points $\lambda\in\sS$ corresponding to a quadrature rule. In fact, if one requires only
matrix-vector multiplies as in Krylov-based iterative methods (e.g.~\textsc{cg}~\cite{Hestenes52} or
\textsc{minres}~\cite{Paige75} solvers) for the Galerkin system of equations, this restriction can be relaxed to
matrix-vector products $A(\lambda)\vv$ for a given $N$-vector $\vv$; there is no need for the coefficients of an
expansion of $A(s)$ as in \eqref{eq:Aexpansion}. Therefore the method takes full advantage of sparsity of the
parameterized system resulting in reduced memory requirements. The decomposition also yields insights for
preconditioning the Galerkin system that generalize existing work. Additionally, the decomposition immediately reveals
bounds on the eigenvalues of the Galerkin system for the case of symmetric $A(s)$.

The remainder of the paper is structured as follows. In Section \ref{sec:galerkin}, we derive the Galerkin method for a
given problem and basis set. We then derive the decomposition of the matrix used to compute the Galerkin coefficients
and examine its consequences including bounds on the eigenvalues of the matrix, strategies for simple implementation
and code reuse, and insights into preconditioning. We then provide a numerical study of various preconditioners
suggested by the decomposition using the common test case of an elliptic PDE with parameterized coefficients in Section
\ref{sec:precon}; the codes for the numerical study are available in a \textsc{Matlab} suite of tools accessible
online~\cite{Gleich10}. To emphasize the advantages of code reuse, we apply the method to an engineering test problem
using existing codes in Section \ref{sec:joe}. Finally, we conclude in Section \ref{sec:summary} with summarizing
remarks.

\subsection{Notation}
For the remainder of the paper, we will use the bracket notation $\ip{\cdot}$ to denote a \emph{discrete} approximation
to the integral with respect to the weight function $\omega$. In other words, for functions
$f:\sS\longrightarrow\mathbb{R}$, 
\begin{equation}
\int_{\sS}f(s)\omega(s)\,ds
\;\approx\;
\sum_{\beta\in\sJ} f(\lambda_\beta)\nu_\beta
\;\equiv\;
\ip{f},
\end{equation}
where the points $\lambda_\beta=(\lambda_{\beta_1},\dots,\lambda_{\beta_d})\in\sS$ and weights
$\nu_\beta\in\mathbb{R}_+$ define a multivariate quadrature rule\footnote{We have restricted our attention to
quadrature rules with positive weights.} for multi-indicies in the set $\sJ\subset\mathbb{N}^d$. We will abuse this
notation by putting matrix and vector valued functions inside the brackets, as well. For example, $\ip{A}$ denotes the
mean of $A(s)$ computed with a quadrature rule.

\section{A Spectral Galerkin Method}
\label{sec:galerkin}

We first review the Galerkin method~\cite{Constantine10c} for computing the coefficients of the polynomial model
\eqref{eq:polymodel}. Define the residual
\begin{equation}
\label{eq:resid}
r(y,s) = A(s)y(s)-b(s),
\end{equation}
and let $x_g(s)$ be the Galerkin approximation. Denote the $i$th component of the residual by $r_i(x_g,s)$.
We require that each component of the residual be orthogonal to the approximation space defined by the span of
$\pi_\alpha$ with $\alpha\in\sI$,
\begin{equation}
\label{eq:orthogresid}
\ip{r_i(x_g)\pi_\alpha} = 0, \qquad i=1,\dots,N,\qquad\alpha\in\sI.
\end{equation}
We can combine the equations in \eqref{eq:orthogresid} using the matrix notation
\begin{equation}
\ip{r(x_g)\bpi^T} \;=\; \ip{(Ax_g-b)\bpi^T} \;=\; \mathbf{0},
\end{equation}
or equivalently, upon substituting the model $x_g(s)=\mX\bpi(s)$,
\begin{equation}
\label{eq:varform}
\ip{A\mX\bpi\bpi^T} = \ip{b\bpi^T}.
\end{equation}
Using the vec notation~\cite[Section 4.5]{GVL96}, we can rewrite
\eqref{eq:varform} as
\begin{equation}
\label{eq:galerkinsys}
\ip{\bpi\bpi^T\otimes A}\vx=\ip{\bpi\otimes b}.
\end{equation} 
where $\vx=\tvec(\mX)$ is an $N|\sI|\times 1$ constant vector equal to the columns of $\mX$ stacked on top of each
other. The constant matrix $\ip{\bpi\bpi^T\otimes A}$ has size $N|\sI|\times N|\sI|$ and a distinct block structure; the
$\alpha,\beta$ block of size $N\times N$ is equal to $\ip{\pi_\alpha\pi_\beta A}$ for multi-indicies
$\alpha,\beta\in\sI$. Similarly, the $\alpha$ block of the $N|\sI|\times 1$ vector $\ip{\bpi\otimes b}$ is equal to
$\ip{b\pi_\alpha}$, i.e.~the Fourier coefficient of $b$ associated with $\pi_\alpha(s)$ approximated with the quadrature
rule.

Much of the literature on PDEs with random inputs~\cite{Pellissetti00,Powell09-2} points out an interesting block
sparsity pattern that arises in the matrix $\ip{\bpi\bpi^T\otimes A}$ when the parameterized matrix $A(s)$ depends at
most linearly on any components of $s$ and integrals are computed exactly; this is a result the
orthogonality of the bases $\bpi(s)$. For general analytic dependence on the parameters, such sparsity patterns do not
appear~\cite{Ernst10}. However, one can always mimick the sparsity pattern of $A(s)$ with a simple reordering of the
variables. By taking the transpose of \eqref{eq:varform} and using the same vec operations, if we define
$\tilde{\mathbf{x}}=\tvec(\mX^T)$ (i.e.~the same unknowns reordered), then
\begin{equation}
\label{eq:galerkinsys2}
\ip{A\otimes\bpi\bpi^T}\tilde{\vx}=\ip{b\otimes\bpi}.
\end{equation} 
Notice that the matrix in \eqref{eq:galerkinsys2} retains the sparsity of $A(s)$ in its blocks, since each
$i,j$ block of size $|\sI|\times|\sI|$ is equal to $\ip{A_{ij}\bpi\bpi^T}$, where $A_{ij}(s)$ is the $i,j$
element of $A(s)$. For the remainder of the paper, however, we will work with the form \eqref{eq:galerkinsys}. 

By writing out the numerical integration rule for the integrals used to form $\ip{\bpi\bpi^T\otimes A}$, we uncover an
interesting decomposition, which we state as a theorem.

\begin{theorem}
\label{thm:decomp}
Let $\{(\lambda_\beta,\nu_\beta)\}$ with $\beta\in\sJ$ be a multivariate quadrature rule. The matrix
$\ip{\bpi\bpi^T\otimes A}$ can be decomposed as
\begin{equation}
\label{eq:decomp3}
\ip{\bpi\bpi^T\otimes A} 
= (\mQ\otimes\mI)A(\blambda)(\mQ\otimes\mI)^T,
\end{equation}
where $\mI$ is the $N\times N$ identity matrix, and $\mQ$ is a matrix of size $|\sI|\times|\sJ|$ -- one row for each
basis polynomial and one column for each point in the quadrature rule. The matrix $A(\blambda)$ is a block diagonal
matrix of size $N|\sJ|\times N|\sJ|$ where each nonzero block is $A(\lambda_\beta)$ for $\beta\in\sJ$.
\end{theorem}

\begin{proof}
Writing out the quadrature rules,
\begin{equation}
\label{eq:decomp1}
\ip{\bpi\bpi^T\otimes A} 
= \sum_{\beta\in\sJ} \left[\bpi(\lambda_\beta)\bpi(\lambda_\beta)^T\otimes A(\lambda_\beta)\right]\nu_\beta.
\end{equation}
Notice that the elements of the vector $\bpi(\lambda_\beta)$ are the polynomials $\pi_\alpha(s)$ evaluated at the
quadrature point $\lambda_\beta$. If we define the vectors
\begin{equation}
\vq_\beta = \sqrt{\nu_\beta} \bpi(\lambda_\beta),
\end{equation}
then we have
\begin{equation}
\label{eq:decomp2}
\ip{\bpi\bpi^T\otimes A} 
= \sum_{\beta\in\sJ} \vq_\beta\vq_\beta^T\otimes A(\lambda_\beta).
\end{equation}
Let $\mQ$ be the matrix whose columns are $\vq_\beta$, and define the block diagonal matrix
$A(\blambda)$ with diagonal blocks equal to $A(\lambda_\beta)$. Then for an $N\times N$ identity matrix $\mI$, we can
rewrite \eqref{eq:decomp2} as \eqref{eq:decomp3}, as required.
\end{proof}

As an aside, we note that if $A(s)$ depends polynomially on the parameters $s$, then each integrand in the matrix
$\ip{\bpi\bpi^T\otimes A}$ is a polynomial in $s$. Therefore by the polynomial exactness, there is a Gaussian
quadrature rule such that the numerical integration approach exactly recovers the true Galerkin matrix. 

The elements of $\mQ$ are the orthogonal polynomials evaluated at the quadrature points and multiplied by the square
root of the quadrature weights. They are intimately related to the normalized eigenvectors of the symmetric,
tridiagonal Jacobi matrices of the three-term recurrence coefficients of the orthogonal polynomials, which can be
computed efficiently by methods for computing eigenvectors; see~\cite{Gautschi04,Constantine10c} and Appendix A for
more details. For two polynomials $\pi_\alpha(s)$ and $\pi_\beta(s)$ with $\alpha,\beta\in\sI$, define $\vr_\alpha$ and
$\vr_\beta$ to be the corresponding rows of $\mQ$; then
\begin{equation}
\vr_\alpha\vr_\beta^T  = \sum_{\gamma\in\sJ} \pi_\alpha(\lambda_\gamma) \pi_\beta(\lambda_\gamma)\nu_\gamma.
\end{equation}
If the quadrature rule is a tensor product Gaussian quadrature rule of sufficiently high order to exactly compute the
polynomial integrand, then this implies $\mQ\mQ^T=\mI$, where $\mI$ is the $|\sI|\times|\sI|$ identity matrix; we will
assume from here on that the chosen quadrature rule yields this property. We will also assume that the number of basis
polynomials is less than the number of points used in the quadrature rule, i.e.~$|\sI|\leq|\sJ|$; otherwise the Galerkin
matrix $\ip{\bpi\bpi^T\otimes A}$ will be rank deficient. 

Theorem \ref{thm:decomp} reveals bounds on the eigenvalues of the matrix $\ip{\bpi\bpi^T\otimes A}$ for the case when
$A(s)$ is symmetric; we state this as a corollary.

\begin{corollary}
\label{cor:eigbounds}
Suppose $A(s)$ is symmetric for all $s\in\sS^d$. The eigenvalues of $\ip{\bpi\bpi^T\otimes A}$ satisfy the bounds
\begin{equation}
\min_{\beta\in\sJ}\left[\theta_{\mathrm{min}}\Bigl(A(\lambda_\beta)\Bigr)\right]
 \;\leq\;
\theta\Bigl(\ip{\bpi\bpi^T\otimes A}\Bigr) 
\;\leq\;
\max_{\beta\in\sJ} \left[\theta_{\mathrm{max}}\Bigl(A(\lambda_\beta)\Bigr)\right],
\end{equation}
where $\theta(X)$ denotes the eigenvalues of a matrix $X$, and $\theta_{\mathrm{min}}(X)$ and $\theta_{\mathrm{max}}(X)$
denote the smallest and largest eigenvalues of $X$, respectively.
\end{corollary}

The proof of Corollary \ref{cor:eigbounds} involves a detailed discussion of the properties of orthogonal polynomials
and the tridiagonal matrices of their three-term recurrence coefficients. To keep this section focused, we have placed
the proof in an appendix along with a note on the sharpness of the bounds. To conclude this section, we note that
if $A(s)$ is positive definite for all $s\in\sS$, then Corollary 2.2 implies that $\ip{\bpi\bpi^T\otimes A}$ will be
positive definite, as well.
% \begin{proof}
% Let $\{(\lambda_\beta,\nu_\beta)\}$ be a tensor product Gaussian quadrature rule of sufficiently high order so that
% $\mQ\mQ^T=\mI$. Let $\tilde{\bpi}(s)$ be a vector containing the polynomial bases such that
% $[\bpi(s)^T,\tilde{\bpi}(s)^T]^T$ contains the tensor product polynomial basis corresponding to the tensor product
% quadrature grid $\{\lambda_\beta\}$. Let $\tilde{\mQ}$ be the matrix with columns $\tilde{\vq}_\beta =
% \tilde{\bpi}(\lambda_\beta)\sqrt{\nu_\beta}$; by construction, the rows of $\tilde{\mQ}$ are orthonormal and
% $\tilde{\mQ}\mQ^T=\mathbf{0}$. Also, by the Christoffel-Darboux formula~\cite[Theorem 1.32]{Gautschi04},
% $\mQ^T\mQ+\tilde{\mQ}^T\tilde{\mQ}=\mI$. Then the matrix
% \begin{equation}
% \mZ = \bmat{\mQ\\ \tilde{\mQ}}
% \end{equation}
% is square and orthogonal, i.e. $\mZ^{-1}=\mZ^T$. Notice that the matrix  
% $(\mQ\otimes\mI)A(\blambda)(\mQ\otimes\mI)^T$
% is the first $|\sI|\times|\sI|$ principal minor of $(\mZ\otimes\mI)A(\blambda)(\mZ\otimes\mI)^T$.
% Then by the interlacing theorem~\cite[Theorem 8.1.7]{GVL96}, the eigenvalues of
% $(\mQ\otimes\mI)A(\blambda)(\mQ\otimes\mI)^T$ are bounded by the extreme eigenvalues of
% $(\mZ\otimes\mI)A(\blambda)(\mZ\otimes\mI)^T$. But since $\mZ^T=\mZ^{-1}$, this is a similarity transformation with
% $A(\blambda)$, which completes the proof.
% \end{proof}

\subsection{Iterative Solvers}

Due to their size and sparsity, the preferred way of solving \eqref{eq:galerkinsys} is with Krylov based iterative
solvers~\cite{Pellissetti00} that rely on matrix-vector products with the matrix $\ip{\bpi\bpi^T\otimes A}$. By
employing the decomposition in Theorem \ref{thm:decomp}, we can compute these using only multiplication of the
parameterized matrix evaluated at the quadrature point by a given vector. More precisely, given a vector
$\vu=\tvec(\mU)$, suppose we wish to compute
\begin{equation}
\vv\;=\;\tvec(\mV)\;=\;(\mQ\otimes\mI)A(\blambda)(\mQ\otimes\mI)^T\vu.
\end{equation}
We accomplish this in three steps using the properties of the Kronecker product:
\begin{enumerate}
  \item \label{step1} $\mW=\mU\mQ$. Let $\vw_\beta$ be a column of $\mW$ with $\beta\in\sJ$.
  \item \label{step2} For each $\beta$, $\vy_\beta=A(\lambda_\beta)\vw_\beta$. Define $\mY$ to be the matrix with  
  columns $\vy_\beta$.
  \item \label{step3} $\mV=\mY\mQ^T$.
\end{enumerate}
Step \ref{step1} can be thought of as pre-processing, and step \ref{step3} as post-processing. In practice, each
row of $\mQ$ may have a Kronecker structure corresponding to the tensor product quadrature rule. In this case, steps
\ref{step1} and \ref{step3} can be computed accurately and efficiently using multiplication methods such
as~\cite{Fernandes98}. The second step requires only constant matrix-vector products where the matrix is $A(s)$
evaluated at the quadrature points. Therefore we can take advantage of a memory-efficient, reusable
interface for the matrix-vector multiplies that will exploit any sparsity in the matrix. We reiterate that this can be
accomplished \emph{without} any knowledge of the specific type of parameter dependence in $A(s)$. 

Each of the three steps individually admits embarrassing parallelization; steps \ref{step1} and \ref{step3} are
matrix-vector multiplies with independent right hand sides, and each $\vy_\beta$ in step \ref{step2} can be computed
independently. However, there is substantial communication necessary between the steps, and this will be the primary
barrier to parallel scalability.

The cost of this method depends on the number of basis polynomials in the approximation and number of points in the
quadrature rule. Steps \ref{step1} and \ref{step3} each require $N|\sI||\sJ|$ multiplies. If a matrix-vector product
with $A(s)$ takes $\mathcal{O}(N)$ operations due to its sparsity pattern, then step \ref{step2} takes
$\mathcal{O}(|\sJ|N)$ operations. Methods based on the series expansion \eqref{eq:Aexpansion} of $A(s)$ can be
significantly cheaper -- depending on the number of terms in \eqref{eq:Aexpansion}, and assuming the so-called triple
products are precomputed before the iterative solver is applied. The trade-off between cost and flexibility must be
assessed per application.

\subsection{Preconditioning Strategies}

The number of iterations required to achieve a given convergence criterion (e.g.~a sufficiently small residual) can be
greatly reduced for Krylov-based iterative methods with a proper preconditioner. In general, preconditioning a system is
highly problem dependent and begs for the artful intuition of the scientist. However, the structure revealed by the
decomposition from Theorem \ref{thm:decomp} offers a number of useful clues.

Suppose we have an $N\times N$ matrix $\mP$ that is easily invertible. We can construct a block-diagonal preconditioner
$\mI\otimes\mP^{-1}$, where $\mI$ is the identity matrix of size $|\sI|\times |\sI|$. If we premultiply the
preconditioner against the factored form of $\ip{\bpi\bpi^T\otimes A}$, we get
\begin{equation} \label{eq:precond}
(\mI\otimes\mP^{-1})(\mQ\otimes\mI)A(\blambda)(\mQ\otimes\mI)^T
=
(\mQ\otimes\mI)(\mI\otimes\mP^{-1})A(\blambda)(\mQ\otimes\mI)^T.
\end{equation}
By the mixed product property and commutativity of the identity matrix, the block-diagonal preconditioner
slips past $\mQ\otimes\mI$ to act directly on the parameterized matrix evaluated at the quadrature points. The blocks on
the diagonal of the inner matrix product are $\mP^{-1}A(\lambda_\beta)$ for $\beta\in\sJ$.  In other words, we can
choose one constant matrix $\mP$ to affect the parameterized system at all quadrature points.

A reasonable and popular choice is the mean $\mP=\ip{A}$; see~\cite{Powell08, Pellissetti00} for a detailed analysis of
this preconditioner for stochastic finite element systems. Notice that this is also the first $N\times N$ block of the
Galerkin matrix. However, if $A(s)$ is very large or has some complicated parametric dependence, then forming the mean
system and inverting it (or computing partial factors) for the preconditioner may be prohibitively expensive. If the
dependence of $A(s)$ on $s$ is close to linear, then $\mP=A\left(\ip{s}\right)$ may be much easier to evaluate and just
as effective.

One goal of the preconditioner is reduce the condition number of the matrix, and one way of achieving this is to reduce
the spread of the eigenvalues. If we knew \emph{a priori} which region of the parameter space produced the extrema of
the parameterized eigenvalues of $A(s)$ (e.g.~the boundaries), then we could choose an appropriate parameter value to
construct an effective preconditioner. Unfortunately, we only get to use one such evaluation. Therefore, if
the largest possible value of the parameterized eigenvalues is very large, we may choose this parameter value.
Alternatively, if the smallest eigenvalue over the parameter space is close to zero (for positive definite systems),
then this may be better reduce the condition number of the Galerkin system than the parameter that produces largest
eigenvalue. In the next section, we explore a few choices for the preconditioner $\mP$ on a standard test problem.

Notice that \eqref{eq:precond} suggests two possible routes for implementing the preconditioner. If only an interface
for a preconditioned matrix-vector multiply is available for the parameterized matrix, then the implementation
corresponding to the right hand side of \eqref{eq:precond} is appropriate. However, this requires $|\sJ|$ applications
of the preconditioner, and in all cases we expect that $|\sI| \leq |\sJ|$. Therefore, our implementation employs the
form on the left hand side of \eqref{eq:precond}, which requires $|\sI|$ applications of the preconditioner and can be
computed efficiently due to its block diagonal structure.

\section{Preconditioning Study}
\label{sec:precon}

\begin{figure}
 \caption{Eigenvalues of the KL expansion.  The numbers indicate the 
 percentage of field energy captured when using that many 
 terms in the expansion.}
 \centering
 \includegraphics[width=0.7\linewidth]{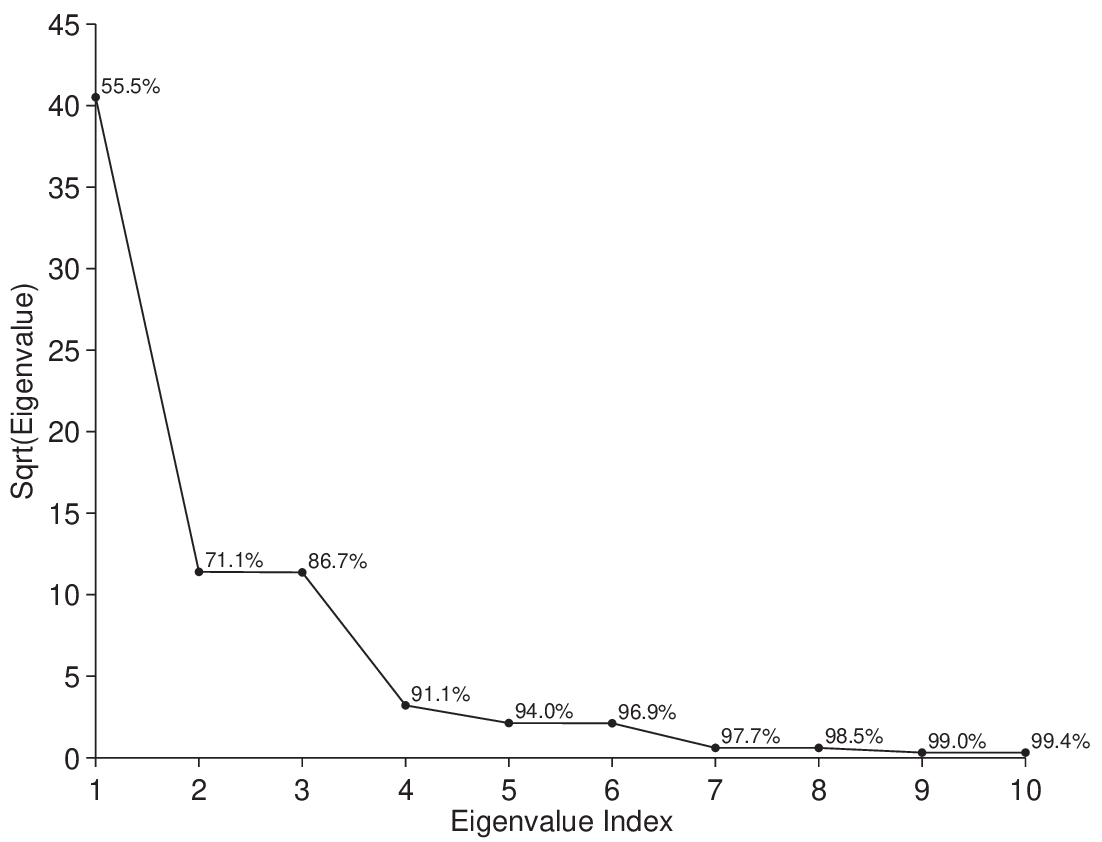}
 \label{fig:kleigs}
\end{figure}

Consider the following parameterized elliptic partial differential equation with homogeneous Dirichlet boundary
conditions; variations of this problem can be found in~\cite{Babuska01,Babuska04,Frauenfelder2005}, amongst others. We
seek a solution $u(x,s)$ that satisfies
\begin{equation}
\label{eq:elliptic}
\nabla\cdot(a(x,s) \nabla u(x,s)) = 1, \qquad x\in[0,1]^2, 
\end{equation}
where $u=0$ on the boundary of the square domain $[0,1]^2$, and the parameter space is the hypercube $s\in[-1,1]^4$
equipped with a uniform measure. The logarithm of the elliptic coefficient is given by a truncated Karhunen-Loeve
like expansion~\cite{Loeve78} of a zero-mean random field with covariance 
\begin{equation}
C(x_1,x_2)=2\exp(-\|x_1-x_2\|^2/2),
\end{equation}
so that
\begin{equation}
\log(a(x,s)) = 2\sum_{k=1}^4 \sigma_k\psi_k(x)s_k,
\end{equation}
where $\{\sigma_k^2,\psi_k\}$ are the eigenpairs of the covariance function.
We discretize \eqref{eq:elliptic} using the finite element method implemented 
in \textsc{Matlab}'s \textsc{pde toolbox} on an
irregular mesh of $N =1,\!921$ triangles. 
We solve the eigenvalue problem with the discrete covariance matrix to compute
the values of the eigenfunctions $\psi_k(x)$ on the given mesh. 
The square roots of first ten eigenvalues $\sigma_k$
are plotted in Figure \ref{fig:kleigs}; we truncate the expansion at $d=4$ to retain 
roughly 90\% of the energy of the field.
Given a point in the parameter space, the PDE Toolbox allows us to 
access the stiffness matrix. 
Therefore we can apply
\textsc{Matlab}'s \textsc{minres} solver using only matrix-vector 
multiplies against the stiffness matrix evaluated at the quadrature
points.
For the polynomial basis, we use the normalized multivariate 
product Legendre polynomials of order 5, which includes
$|\sI|=\binom{4+5}{5}=126$ basis functions. Therefore the number of unknowns 
is $N|\sI|=242,\!046$. We use a tensor product
Gauss-Legendre quadrature rule of order 12 in each parameter 
(20,736 points) to implicitly form the Galerkin matrix;
this is more than sufficient to maintain the 
orthogonality in the basis functions. 

\begin{figure}
 \centering
 \includegraphics[width=0.8\linewidth]{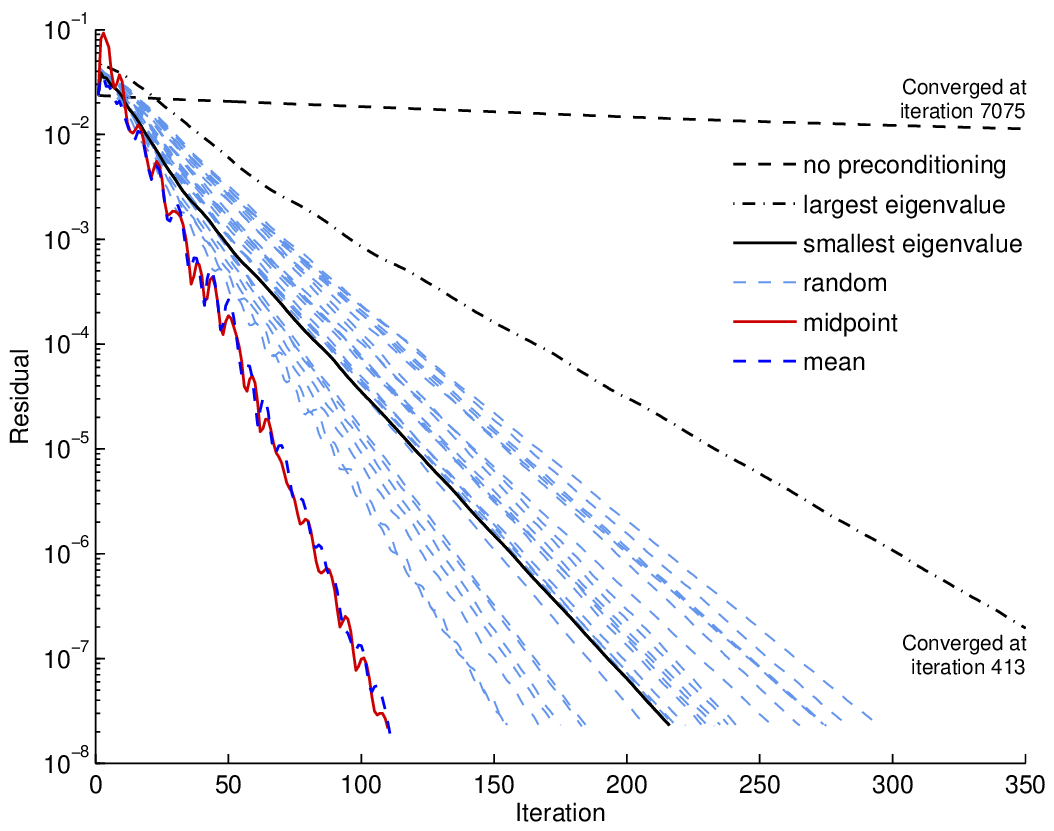}
 \caption{Convergence of the 2-norm of the residual for a variety of preconditioners.}
 \label{fig:precond}
\end{figure}

To test the preconditioning ideas suggested by the decomposition from Theorem \ref{thm:decomp}
and equation~\eqref{eq:precond}, we try five different choices of $\mP$. For each $\mP$, we precompute the Cholesky
factors to apply \eqref{eq:precond} efficiently. Figure~\ref{fig:precond} and Table~\ref{tab:precond-time} summarize
the results, which we now explain in detail.  
In the first type of preconditioning, we set $\mP = A(s_r)$ for a random point $s_r$ in $[-1,1]^4$.  
To sample the likely effectiveness of any random point, we select 25 random points and
evaluate each choice of $s_r$.  
Second, we set $\mP=A(s_{\max})$ where $A(s_{\max})$ is the matrix with the largest eigenvalue
in $[-1,1]^4$.  To compute $s_{\max}$, we use up to 100 iterations of the power method
to estimate the largest eigenvalue at each point in the quadrature rule, that is
for each $\lambda_\beta$ for $\beta \in \mathcal{J}$;  we also evaluate $A(s)$ for
the tensor product of the endpoints as well.  Third, we set $\mP=A(s_{\min})$ where
$A(s_{\min})$ is the matrix with the smallest eigenvalue in $[-1,1]^4$.  For
this computation, we use \textsc{Matlab}'s optimization routine \verb#fmincon#
and use \verb#eigs#/\textsc{arpack} to evaluate the smallest eigenvalue \cite{Lehoucq98}.  
Fourth, we set $\mP=A(s_\text{mid})$ where $s_\text{mid}$ is the midpoint of the domain,
which is the origin for our experiments.  Fifth, and finally, we set $\mP=\ip{A}$,
the mean preconditioner.  To evaluate the mean, we use a 2nd order quadrature rule for a 
fast approximation and a 5th order quadrature rule for a more exact approximation.

In Figure~\ref{fig:precond}, we show the convergence of the $2$-norm of the residual
for each of the preconditioning strategies, as well as no preconditioning.  We only
show the result from the 5th order mean based preconditioner as there was
no appreciable difference in convergence.  
This plot clearly shows that the mean and midpoint preconditioners are excellent choices.
Further note that any random point is more than two orders of magnitude better than
no preconditioning at all.  In fact, using a random point is better than using
the point with largest eigenvalue and may compare with using the point with 
the smallest eigenvalue.

Next, Table~\ref{tab:precond-time} shows a few timing results from these
experiments.  We present the time taken to compute/setup the preconditioner,
the average time between iterations, and the total time taken by the \textsc{minres}
method (excluding preconditioner setup).
The times are from \textsc{Matlab} 2010a on an Intel Core i7-960
processor (3.2 GHz, 4 cores) with 24GB of RAM.  Matlab used
four cores for its own multithreading and we used the Parallel Computing
toolbox's \verb#parfor# construction to parallelize the matrix-vector product
over 4 cores as well.  By watching a processor performance meter, we observed high
utilization of all four cores.  However, the system's memory was nearly exhausted by the experiment.
At irregular intervals, the system would begin swapping memory
to disk heavily. This behavior caused erratic timing results, especially
for the two random point evaluations. On repeated runs of the experiments,
we observed wall-clock time differences of up to 100 seconds
for any individual result. Thus differences below this threshold are
not meaningful. From these timing results,
we conclude that preconditioning changes the iteration time only
slightly -- if at all.  Furthermore, the setup times are often small when compared
with the savings in runtime.  

\begin{table}
 \centering
 \caption{Timing results for preconditioning showing only two of the 
 random points.}
 \label{tab:precond-time}
 \begin{tabularx}{\linewidth}{lXXXX}
\toprule
Method & Setup (s) & Iterations & Avg.~Iter.~(s) & Total (s)\\
\midrule
  no preconditioning &    0.0 &  7075 &    6.5 & 46251.1 \\
  largest eigenvalue &  176.0 &   413 &    6.6 &  2810.7 \\
 smallest eigenvalue &   16.0 &   216 &    6.5 &  1469.4 \\
              random &    0.1 &   238 &    6.7 &  1661.8 \\
              random &    0.9 &   216 &    8.4 &  2030.2 \\
    mean (5th order) &    4.7 &   111 &    6.7 &   807.2 \\
    mean (2nd order) &    0.7 &   110 &    6.8 &   820.4 \\
            midpoint &    0.1 &   110 &    7.0 &   843.2 \\
\bottomrule
 \end{tabularx}

\end{table}

The Matlab codes for this study include a utility for generating realizations of random fields~\cite{Constantine10-2}
and a suite of tools called PMPack (Parameterized Matrix Package) for using spectral methods to approximate the solution
of parameterized matrix equations.  The codes for generating the results of the study can be found at~\cite{Gleich10}.

\section{Application -- Heat Transfer with Uncertain Material Properties}
\label{sec:joe}

As a proof of concept, we examine an application from computational fluid dynamics with uncertain model inputs. The
flow solver used to compute the deterministic version of this problem -- i.e. for a single realization of the model
inputs -- was developed at Stanford's Center for Turbulence Research as part of the Department of Energy's Predictive
Science Academic Alliance Program; the numerical method used is described in~\cite{Pecnik08} and is based on an
implicit, second order spatial discretization. For this example, we slightly modified the codes to extract of the
non-zero elements of the matrix and right hand side used in the computation of the temperature distribution. With
access to the matrix-vector multiply, we were able to apply the Galerkin method to the stochastic version of the
problem to approximate the statistics of solution.

\subsection{Problem Set-up}

\begin{figure}
\begin{center}
\includegraphics[scale=0.5]{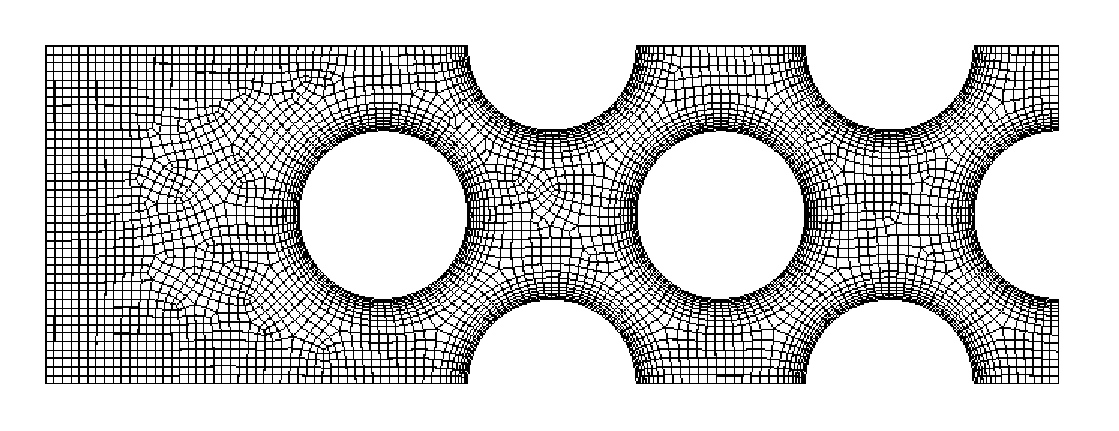}
\end{center}
\caption{Mesh used to discretize the domain $\Omega$ and compute temperature distribution.}
\label{fig:joedomain}
\end{figure}

The governing equation is the integral version of a two-dimensional steady advection-diffusion equation. We
seek a scalar field $\phi=\phi(x,y)$ representing the temperature defined on the domain $\Omega$ that satisfies
\begin{equation}
\int_{\partial\Omega} \rho\phi\left(\vec{v}(s)\cdot\vec{dS}\right)
=
\int_{\partial\Omega}
\left(\Gamma(s)+\Gamma_t\right)\left(\nabla\phi\cdot\vec{dS}\right).
\end{equation}
The density $\rho$ is assumed constant. The velocity $\vec{v}$ is precomputed by solving the incompressible
Reynolds averaged Navier-Stokes equations and randomly perturbed by three spatially varying oscillatory functions with
different frequencies such that the divergence free constraint is satisfied; the magnitudes of the perturbations are
parameterized by $s_1$, $s_2$, and $s_3$, respectively. We interpret the magnitudes as uniform random perturbations of
the velocity field, which is simply an input to this model. The diffusion coefficient $\Gamma=\Gamma(s_4,s_5,s_6)$ is
similarly altered by a strictly positive, parameterized, spatially varying function that models random perturbation.
Collectively, the parameters $s_1,\dots,s_6$ are independent and distributed uniformly over $[-1,1]$; the
parameter space becomes the hypercube $[-1,1]^6$. The turbulent diffusion coefficient $\Gamma_t$ is evaluated according
to the Spalart-Allmarass model~\cite{Pecnik08}. The domain $\Omega$ is a channel with a series of cylinders; the
computational mesh on the domain $\Omega$ contains roughly 10,000 nodes and is shown in Figure \ref{fig:joedomain}.
Inflow and outflow boundary conditions are prescribed in the streamwise direction, and periodic boundary conditions are
set along the $y$ coordinate. Specified heat flux boundary conditions are applied on the boundaries of the cylinders to
model a cooling system.

The goal is to compute the expectation and variance of the scalar field $\phi=\phi(x,y,s)$ over the domain
$\Omega$ with respect to the variability introduced by the parameters. We use the Galerkin method to construct a
polynomial approximation of $\phi$ along the coordinates induced by the parameters $s$ using product Legendre
polynomials. To show that this method applies for an arbitrary basis set, we choose the basis polynomials to reflect the
solution's anisotropic dependence on the parameters; it is neither the standard full polynomial
or tensor product polynomial basis. The order of univariate polynomial associated with each parameter is given in Table
\ref{tab:polyorder}, and multivariate bases are included to fall within a convex index set. The basis contains 142
multivariate polynomials; for a detailed description of the choice of bases, see~\cite{Constantine09-3}. To solve the
Galerkin system, we use \textsc{Matlab}'s \textsc{BiCGstab}~\cite{Vorst1992} method (since the matrix is not
symmetric). For a preconditioner, we could only access the diagonal elements of the matrix from the solver. The results
of the preconditioning study in Section \ref{sec:precon} encouraged us to choose the midpoint of the hypercube
parameter space to construct the preconditioner.

\begin{table}
\begin{center}
\begin{tabular}{cccccc}
$s_1$ &$s_2$ &$s_3$ &$s_4$ &$s_5$ &$s_6$ \\
\hline
3 & 1 & 1 & 8 & 5 & 5
\end{tabular}
\end{center}
\caption{The order of univariate polynomial for each variable used in the basis set.}
\label{tab:polyorder}
\end{table}

\begin{figure}
\begin{center}
\includegraphics[scale=0.4]{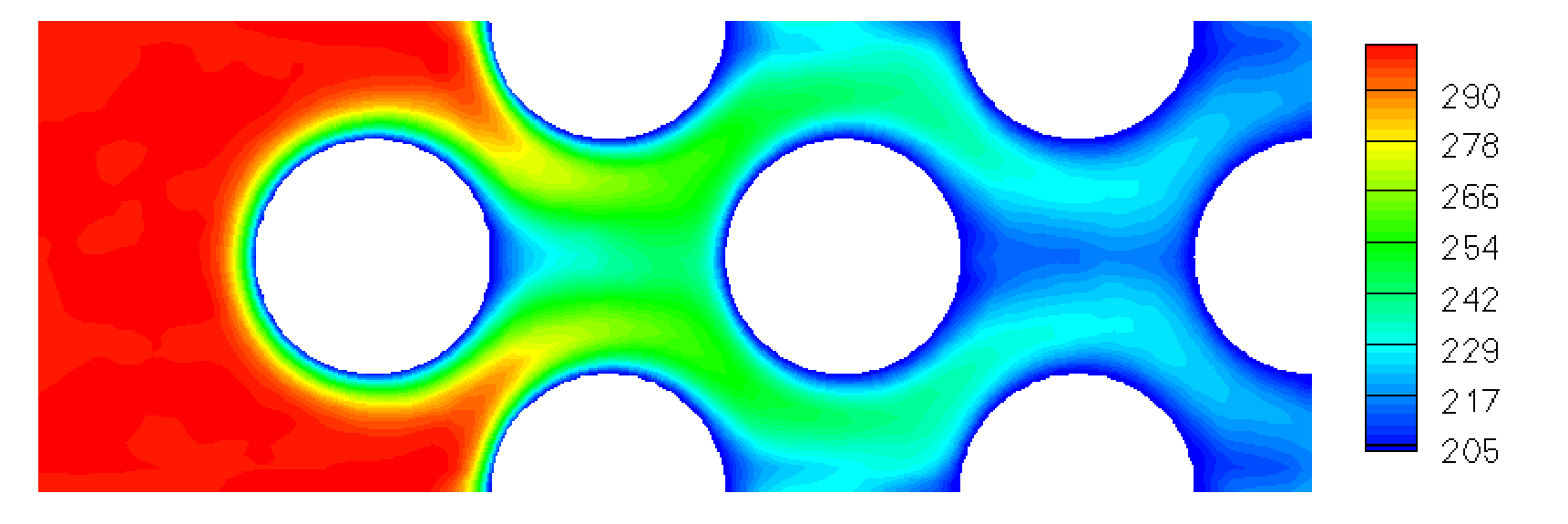}
\includegraphics[scale=0.4]{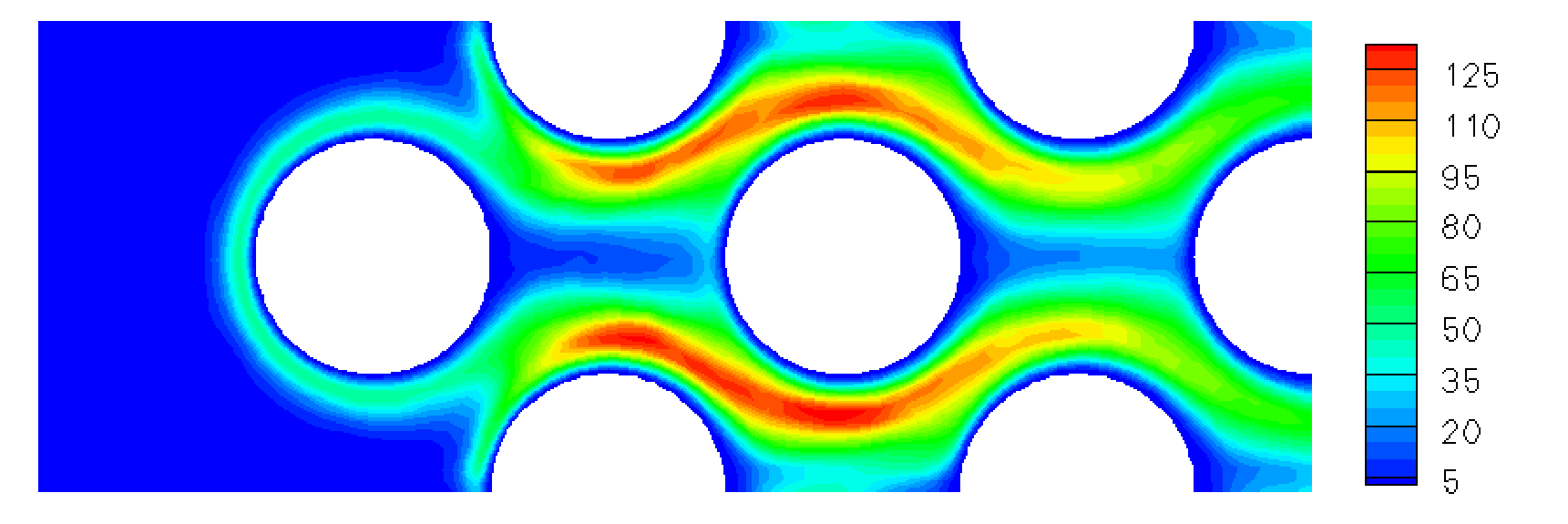}
\end{center}
\caption{The expectation (above) and variance (below) of the temperature field
$\phi$ over the domain. Red corresponds to larger values and blue corresponds to smaller values.}
\label{fig:joestats}
\end{figure}

We plot the expectation and variance of $\phi$ over the domain $\Omega$ in Figure \ref{fig:joestats}, which are computed
in the standard way as functions of the Galerkin coefficients. The variance in $\phi$ occurs in the
downstream portion of the domain as a result of the variability in the diffusion coefficient.

\section{Summary}
\label{sec:summary}

We have examined the system of equations arising from a spectral Galerkin approximation of the vector valued solution
$x(s)$ to the parameterized matrix equation $A(s)x(s)=b(s)$. Such problems appear when PDE models with parameterized (or
random) inputs are discretized in space, and a Galerkin projection with an orthogonal polynomial basis is used for
approximation in the parameter space. We showed how the system of equations used to compute the coefficients of the
Galerkin approximation admits a factorization once the integration is formally replaced by a numerical quadrature rule
-- a common step in practice. The factorization involves (i) a matrix with orthogonal rows related to the chosen
polynomial basis and numerical quadrature rule and (ii) a block diagonal matrix with nonzero blocks equal to $A(s)$
evaluated at the quadrature points. Then matrix-vector products with the Galerkin system can be computed with only the action of $A(s)$
on a vector at a point in the parameter space; this yields a reusable interface for implementing the Galerkin method.
The factorization also reveals bounds on the spectrum of the Galerkin matrix, and the Kronecker structure of the
factorization gives clues to successful preconditioners. We tested some ideas from these preconditioner clues on the
standard test problem of an elliptic PDE with parameterized coefficients, and we saw that using $A(s)$ evaluated at the
midpoint of the parameter space was as effective as the popular mean-based preconditioner; the midpoint preconditioner
is much easier to compute, in general. As a proof of concept, we applied the method to an engineering flow problem by
slightly modifying an existing CFD code to retreive the matrix elements.

\section{Acknowledgments}
This material is based upon work supported by the Department of Energy [National Nuclear Security Administration] under
Award Number NA28614. Sandia National Laboratories is a multi-program laboratory operated by Sandia Corporation, a
wholly owned subsidiary of Lockheed Martin company, for the U.S. Department of Energy’s National Nuclear Security
Administration under contract DE-AC04-94AL85000.

\section*{Appendix A: Proof of Corollary \ref{cor:eigbounds}}

To prove Corollary \ref{cor:eigbounds}, we proceed somewhat circuitously as follows. Much of the notation for the
orthogonal polynomials in this proof is taken from~\cite{Gautschi04}. Throughout the proof, we will use the index $i$ to
refer to a specific parameter, so that $i=1,\dots,d$. 

Let $\mJ_i$ be the $n_i \times n_i$ symmetric, tridiagonal Jacobi matrix of three-term recurrence coefficients for the
univariate polynomials $\bpi_i(s_i)$ which are orthogonal with respect to the weight function $\omega_i(s_i)$; the
vector $\bpi_i(s_i)$ contains the polynomials up to order $n_i-1$ arranged in ascending degree from top to bottom. The
three-term recurrence relation for the polynomials can be written in matrix form as
\begin{equation}
s_i\bpi_i(s_i) = \mJ_i \bpi_i(s_i) + \beta \pi_{n_i}(s_i)\ve, 
\end{equation}
where $\pi_{n_i}(s_i)$ is the univariate orthogonal polynomial of order $n_i$, $\beta$ is a constant that completes the
three-term recurrence relationship, and $\ve$ is an $n_i$-vector of zeros with a one in the last entry. Let $\lambda$ be
a zero of $\pi_{n_i}(s_i)$, so that
\begin{equation}
\lambda\bpi_i(\lambda) = \mJ_i \bpi_i(\lambda).
\end{equation}
This immediately yields an eigenpair $\{\lambda,\bpi_i(\lambda)\}$ for $\mJ_i$. Since there are $n_i$ zeros of the
univariate polynomial $\pi_{n_i}(s_i)$, we have a complete set of eigenpairs for $\mJ_i$, and note that the eigenvalues
of $\mJ_i$ are the points of the $n_i$-point Gaussian quadrature rule for the measure $\omega_i(s_i)$. The weights of
the rule are given by the square roots of the first elements of the normalized eigenvector, which is
$1/\|\bpi_i(\lambda_k)\|$ for $k=1,\dots,n_i$. Let $\mZ_i$ be the matrix whose $k$th column is given (in MATLAB
notation) by
\begin{equation}
\mZ_i(:,k) = \frac{1}{\|\bpi_i(\lambda_k)\|}\bpi_i(\lambda_k),
\end{equation}
so that $\mZ_i$ contains the normalized eigenvectors of $\mJ_i$, i.e. $\mZ_i^T=\mZ_i^{-1}$. Next define 
\begin{equation}
\mZ = \mZ_1\otimes\cdots\otimes\mZ_d,
\end{equation}
where $\otimes$ denotes the Kronecker product of matrices. By the mixed product property, $\mZ^T=\mZ^{-1}$. Notice that
the elements of $\mZ$ can be referenced by a multi-index, i.e.
\begin{equation}
\mZ_{\alpha\beta} = \sqrt{\nu_{\beta}}\pi_{\alpha_1}(\lambda_\beta)\cdots\pi_{\alpha_d}(\lambda_\beta),\qquad
0\leq\alpha_i \leq n_i-1, \qquad i=1,\dots,d,
\end{equation}
and $\{(\lambda_\beta,\nu_\beta)\}$ with $\beta\in\sJ$ are the point/weight pairs of a tensor product Gaussian
quadrature rule of order $n_i$ in variable $i$. Let $\tilde{\bpi}(s)$ be the multivariate product orthogonal
polynomials corresponding with the rows of $\mZ$. In other words, $\tilde{\bpi}(s)$ contains the multivariate
orthogonal polynomials corresponding to the index set
\begin{equation}
\tilde{\sI} = \{\alpha \: | \: \max_i \alpha_i < n_i,\; i=1,\dots,d\}.
\end{equation}
Assume without loss of generality that the polynomials $\bpi(s)$ used in the spectral Galerkin method are a subset of
$\tilde{\bpi}(s)$ so that $\sI\subseteq\tilde{\sI}$. Then also
\begin{equation}
\mQ = \Pi\mZ,
\end{equation}
where $\Pi$ is a $|\sI|\times|\sJ|$ selector matrix of zeros and ones. Using Theorem \ref{thm:decomp},
\begin{align*}
\ip{\bpi\bpi^T\otimes A}  &= (\mQ\otimes\mI) A(\blambda) (\mQ\otimes\mI)^T\\
&= (\Pi\otimes \mI)(\mZ\otimes\mI) A(\blambda) (\mZ\otimes\mI)^T (\Pi\otimes\mI)^T.
\end{align*}
Therefore, $\ip{\bpi\bpi^T\otimes A}$ is a principal minor of the matrix $(\mZ\otimes\mI) A(\blambda)
(\mZ\otimes\mI)^T$, which is a similarity transformation of $A(\blambda)$. Since $A(s)$ is symmetric, an interlacing
theorem~\cite[Theorem 8.1.7]{GVL96} tells us that the eigenvalues of $\ip{\bpi\bpi^T\otimes A}$ are bounded by the
extreme eigenvalues of $A(\blambda)$. Finally note that these bounds are sharp when $\sI=\tilde{\sI}$, i.e. when the
basis polynomials in the Galerkin approximation are equal to the tensor product polynomials corresponding to the tensor
product Gaussian quadrature rule.

\bibliographystyle{siam}
\bibliography{paulconstantine}

\end{document}